\documentclass[10pt]{article}
\usepackage{graphicx}
\usepackage{psfrag}
\usepackage{epsf}
\usepackage[cp1250]{inputenc}
\usepackage{amsmath,amsfonts,amssymb,latexsym}
\usepackage[left=2cm,top=2cm,right=2cm,bottom=4cm,nohead]{geometry}
\usepackage{setspace}
\usepackage{enumitem}
\usepackage{multirow}

\newtheorem{claim}{Claim}

\def\vt{t\kern-0.22em\raise.18ex\hbox{\char'47}\lower.18ex\hbox{}\kern-0.08em}
\newcommand{\old}[1]{{}}

\title{Threshold functions for distinct parts: revisiting Erd\H os--Lehner}

\author{\'Eva Czabarka \\
 Department of Mathematics, University of South Carolina\\
Columbia, SC 29208, USA\\
{\tt czabarka@math.sc.edu}
\and
Matteo Marsili\\
The Abdus Salam International Centre for Theoretical Physics\\ Strada Costiera 11, 34014, Trieste, Italy\\
{\tt marsili@ittp.it}
\and
L\'aszl\'o A. Sz\'ekely\thanks{The first and the  third 
authors acknowledge financial support from the grant \#FA9550-12-1-0405 from the U.S.
Air Force Office of Scientific Research (AFOSR) and the Defense Advanced
Research Projects Agency (DARPA). The third author was also  supported in part by the NSF DMS contract 1000475.
This research started at the MAPCON 12 conference,
where the authors enjoyed the hospitality of the Max Planck Institute for the Physics of Complex Systems, Dresden, Germany. We thank Danny Rorabaugh for his comments.}\\
 Department of Mathematics, University of South Carolina\\
Columbia, SC 29208, USA\\
{\tt szekely@math.sc.edu}
 }


\date{\today}

\newcommand{\eopf}{\raisebox{0.8ex}{\framebox{}}}

{\noindent{\bf Proof.}\ }%
{\hfill\eopf\par\bigskip}%

{\noindent{\bf Proof of #1.}\ }%
{\hfill\eopf\par\bigskip}%

{\noindent{\bf Sketch of a proof.}\ }%
{\hfill\eopf\par\bigskip}%

{\noindent{\bf Proof in progress.}\ }%
{\hfill\eopf\par\bigskip}%
\begin{document}

\maketitle

\begin{abstract}
We study four problems: put  $n$ distinguishable/non-distinguishable balls into $k$ non-empty distinguishable/non-distinguishable
boxes randomly. What is the threshold function $k=k(n) $ to make almost sure that no two boxes contain the same number of balls? The non-distinguishable ball
problems are very close to the Erd\H os--Lehner asymptotic formula for the number of partitions of the integer $n$ into $k$ parts with $k=o(n^{1/3})$. The problem is motivated by the statistics of an experiment, where we only can tell whether
outcomes are identical or different.  
\end{abstract}

\eject

\section{Motivation}
Consider a generic experiment where the state of a complex system is probed with repeated experiments providing 
outcome sequence $x_1,x_2,\ldots,x_n$. 
The experimenter  can tell only whether two outcomes are identical or different. 
So each outcome can be thought of as a sample of i.i.d. draws from an unknown probability distribution 
over a discrete space.  The order of outcomes carry no valuable information for us.
We want to understand this process under the condition that $k$ different outcomes are observed out of $n$ 
experiments.
As an example, think of a botanist collecting specimen of flowers in a yet unexplored forest. In order to 
classify observations into species, all that is needed is an objective criterion to decide whether two 
specimens belong to the same species or not. More generally, think of unsupervised data clustering of a series of observations. 

In the extreme situation where $k=n$ and all outcomes are observed only once, the classification is not very informative. 
At the other extreme, when all the outcomes belong to the same class, the experimenter will think he/she 
has discovered some interesting regularity. In general, the information that a set of $n$ repeated experiments yields, 
is the size $m_i$ of each class (or cluster), i.e. the numbers $m_i:=\vert\{\ell:x_\ell=i\}\vert$ (with $\sum_{i=1}^k m_i=n$), as the relative frequency of the observations is what allows to make comparative statements. 
We expect that when the number $k$ of classes is large there will be several classes of the same size, i.e. that
will not be discriminated by the experiment, whereas when $k$ is small each class will have a different 
size$^*$\footnote{$^*$This observation can be made precise in information theoretic terms. 
The label $X$ of one outcome, taken at random from the sample, is a random variable whose 
entropy $H[X]$ quantifies its information content. The size $m_X$ of the class containing $X$, 
clearly has a smaller entropy $H[m]\le H[X]$, by the data processing inequality \cite{Cover}. When $k$ is small,
we expect that $H[X]=H[m]$ whereas when $k$ is large, $H[X]>H[m]$.}.

This is clearly a problem that can be rephrased in terms of distributions of balls (outcomes) into 
boxes (classes). We consider, in particular, the null hypothesis of random placement of balls into 
boxes. In this framework, the question we ask is what is the critical number of boxes $k_c(n)$ 
such that for $k\ll k_c$ we expect  to find that all boxes $i$ contain a different number $m_i$ of balls 
whereas for $k\gg k_c$ boxes with the same number of balls will exist with high probability.

\section{Introduction}
Recall the surjective version (no boxes are empty) of the twelvefold way of counting  \cite{stanley} p.41: putting $n$ distinguish-able/non-distinguishable balls into $k$ non-empty distinguishable/non-distinguishable
boxes correspond to four basic problems in combinatorial enumeration according to Table~1. 
Our concern in all four type of problems
is the threshold function  $n=n(k) $ that makes almost sure that no two boxes contain the same number of balls for a randomly and uniformly selected ball placement. 
Although studying all distinct parts is a topical issue for integer partitions \cite{andrews} and compositions \cite{hitzenko}, 
it is hard to find  any corresponding results for surjections and set partitions except \cite{elect}. 
Our results regarding the threshold functions are summed up in Table~1 in parentheses. We have to investigate only three problems,
since every $k$-partition of an $n$-element set corresponds to exactly $k!$ surjections from $[n]$ to $[k]$---namely those surjections, whose
inverse image partition is the $k$-partition in question. Also, the threshold function is the same for compositions and partitions, although the number of compositions
corresponding to a partition may vary from 1 to $k!$ Therefore the threshold function for set partitions is the same as the threshold function for surjections.

Our proofs use the first and  the second moment method, the second moment method in the  form (due to Chung and  Erd\H os, see  \cite{span} p.76) below.
 For events $A_1,A_2,...,A_N$, the following inequality holds:
\begin{equation} \label{also}
P\Biggl(\bigcup_{i=1}^N A_i\Biggl)\geq   \frac{ \Bigl( \sum_{i=1}^N P(A_i) \Bigl)^2}{ \sum_{i=1}^N P(A_i)+ 2\sum_{1\leq i<j\leq N } P(A_i\cap A_j)}.
\end{equation}
  We always will assume 
  \begin{equation} \label{natural}
   n\geq \binom{k+1}{2}>\frac{k^2}{2},
  \end{equation} otherwise the parts clearly cannot have distinct sizes. We did not even attempt to
obtain a limiting distribution $f(c)$ when $n$ is $c$ times the threshold function---though such an estimate should be possible to obtain.
Although there are deep asymptotic results on random functions and set partitions based on generating functions (e.g. see Sachkov \cite{sachkov}), we do not see how to apply generating functions for our threshold problems. Our results corroborate some formulae of Knessl and Kessler \cite{knessl}, who used techniques from 
applied mathematics to obtain heuristics for partition asymptotics from basic partition recursions. 
\begin{table} \label{first}
\begin{center}
\begin{tabular}{|cc|c|c|}
\cline{1-4}
& &\multicolumn{2}{|c|}{$k$ non-empty boxes}\\
\cline{3-4}
& & distinguishable& non-distinguishable\\
\cline{1-4}
\multicolumn{1}{|c|}{\multirow{2}{*}{$n$ balls}}& \multicolumn{1}{|c|}{distinguishable}& surjections ($n=k^5$)& set partitions ($n=k^5$)\\
\cline{2-4}
& \multicolumn{1}{|c|}{non-distinguishable}& integer compositions ($n=k^3$) & integer partitions ($n=k^3$)\\
\cline{1-4}
\end{tabular}
\end{center}
\caption{Threshold functions for distinct parts for the four  surjective cases in the twelvefold way of counting.}
\end{table}

\section{Threshold function for integer compositions} \label{compmom}
Recall that the number of compositions of the integer $n$ into $k$ positive parts is ${\cal C}(n,k)={n-1\choose k-1}$.  
Let $A_{ij}(t)$ denote the event 
that the $i^{th}$ and $j^{th}$ parts are equal $t$ in a random composition of $n$ into $k$ positive parts. 
Using the first moment method, it is easy to see that
\begin{eqnarray*}
P(\exists \hbox{ equal parts})& =& P\Biggl( \bigcup_{i<j}\bigcup_t A_{ij}(t) \Biggl) \leq \sum_{i<j}\sum_t P(A_{ij}(t))=
 {k\choose 2} \sum_{t\geq 1} \frac{{\cal C}(n-2t,k-2) }{ {\cal C}(n,k) }\\
&=&
{k\choose 2} \sum_{t\geq 1} \frac{\binom{n-2t-1}{k-3}} {\binom{n-1}{k-1}}
 \leq  {k\choose 2} \frac{{n-2\choose k-2}}{{n-1\choose k-1}}=o(1),
\end{eqnarray*}
as $n/k^3\rightarrow \infty$. We make an elementary claim here that we use several times and leave its proof to the Reader.
\begin{claim} \label{elementary} Assume that we have an infinite list of finite sequences  of non-negative numbers, $a_1(n), a_2(n),...,a_{N(n)}(n)$ for $n=1,2,...$, such that
none of the sequences is identically zero. Assume that the number of increasing and decreasing intervals of these finite sequences is bounded, and that
$\max_i a_i(n)=o\Bigl( \sum_i a_i(n)\Bigl)$ as $n\rightarrow\infty$. Then, for any fized $k$,   as $n\rightarrow\infty$.  we have 
$$\sum_{i: k|i} a_i(n)=\Bigl(\frac{1}{k}+o(1)\Bigl) \sum_i a_i(n).$$
\end{claim}
Using the claim we can get a more precise estimate
\begin{equation} \label{pontosabb}
 \sum_{t\geq 1} {\binom{n-2t-1}{k-3}}= \frac{1+o(1)}{2} {n-2\choose k-2}.
\end{equation}
Next we use (\ref{also}) and (\ref{pontosabb}) to show that $P(\exists \hbox{ equal parts})\rightarrow 1$
as $n/k^3\rightarrow 0$. The numerator of (\ref{also}) is the square of
\begin{equation} \label{num-comp}
\sum_{i<j}\sum_t P(A_{ij}(t))= {k\choose 2} \sum_{t\geq 1} \frac{\binom{n-2t-1}{k-3}} {\binom{n-1}{k-1}}=(1+o(1))\frac{k^3}{4n}
\end{equation}
that grows to infinity. Therefore we can neglect the same term without square in the denominator of (\ref{also}). A second negligible term
in the denominator arises if $i<j$ and $u<v$ make only 3 distinct indices (note that they must occur with the same $t$). The corresponding sum of the probabilities is estimated by
\begin{equation} \label{only3-comp}
 k{k-1\choose 2}\sum_{t\geq 1}    \frac{{\cal C}  (n-3t,k-3)}{   {\cal C}  (n,k)}\leq   \frac{k^3}{ 2}\sum_{t\geq 1}    
 \frac{\binom{n-3t-1}{k-4}}{\binom{n-1}{k-1}}<\frac{k^3}{ 2} \frac{\binom{n-3}{k-3}}{\binom{n-1}{k-1}}<\frac{k^5}{n^2}=o\Biggl(\Bigl(\frac{k^3}{n}\Bigl)^2
 \Biggl).
\end{equation}
A third negligible term
in the denominator arises if all  four indices are distinct, but the four corresponding parts are all the same:
$$\binom{k}{4}
\sum_{t\geq 1} \frac{  {\cal C}  (n-4t,k-4)}{  {\cal C}  (n,k)}.$$
This can be estimated by $O(k^7/n^3)$ like the estimate in
(\ref{only3-comp}), and is similarly negligible. 
The significant term in the denominator is 
\begin{equation} \label{4-comp}
 {k\choose 2}{k-2\choose 2}\sum_{\ell\geq 1}\sum_{t\geq 1,t\not=\ell}    \frac{{\cal C}  (n-2\ell-2t,k-4)}{   {\cal C}  (n,k)}
 \end{equation}
 corresponding to the cases analysis when 2-2 parts are the same. This term will not change asymptotically
when we add the $t=\ell$ cases to the summation. So (\ref{4-comp}) is asymptotically equal to (using Claim~\ref{elementary} again)
\begin{equation} \label{mainterm-comp}
 \frac{k^4}{4}\sum_{\ell\geq 1}\sum_{t\geq 1}    \frac{\binom{n-2\ell-2t-1}{k-5}}{\binom{n-1}{k-1}}\sim \frac{k^4}{8}\sum_{\ell\geq 1}
 \frac{\binom{n-2\ell-2}{k-4}}{\binom{n-1}{k-1}}\sim \frac{k^4}{16}\frac{\binom{n-3}{k-3}}{\binom{n-1}{k-1}}\sim \frac{k^6}{16n^2}.
  \end{equation}

  
  We conclude  that the  numerator in (\ref{also}) is asymptotically equal to its denominator (\ref{mainterm-comp}), proving that\\
  $P(\exists \hbox{ equal parts})\rightarrow 1$
as $n/k^3\rightarrow 0$. 

\section{Erd\H os--Lehner and the  threshold function for integer partitions} 
Let ${\cal D}(n,k) $ denote the number of compositions of $n$ into $k$ distinct positive terms. In the previous section we proved 
\begin{eqnarray} \label{lenn}
\lim_{n/k^3\rightarrow 0} \frac{{\cal D} (n,k)}{{\cal C} (n,k)}&=&1 \hbox{ and} \\ 
\lim_{n/k^3\rightarrow \infty} \frac{{\cal D} (n,k)}{{\cal C} (n,k)}&=& 0.\label{fenn}
\end{eqnarray}
Let $p(n,k)$ denote the number of partitions of $n$ into $k$ positive terms and $q(n,k)$ denote the number of partitions of $n$ into $k$ distinct positive terms. 
For $1\leq x_1\leq x_2\leq ... \leq x_k$,
  the well-known bijection  $x_1+ x_2+ ... + x_k\rightarrow  (x_1)+(x_2+1)+(x_3+2)+...+(x_k+k-1)$ shows that  $q(n,k)=p\bigl(n-\binom{k}{2},k\bigl)$.  

A theorem of  Erd\H os and Lehner (\cite{erdoslehner}, see also in \cite{andrews}) asserts that 
for $k=o(n^{1/3})$, the following asymptotic formula holds:
\begin{equation} \label{EL}
p(n,k)\sim \frac{1}{k!}{{\cal C}(n,k)}=\frac{1}{k!}\binom{n-1}{k-1}.
\end{equation}
Gupta's proof to Erd\H os--Lehner (\cite{gupta}, see also in \cite{andrews}) 
obtains
\begin{equation}\label{guptaproof}
\frac{1}{k!}{{\cal C}(n,k)}\leq p(n,k)=q\Bigl( n+\binom{k}{2},k    \Bigl)\leq \frac{1}{k!}{{\cal C}\Bigl(n+\binom{k}{2},k\Bigl)}
\end{equation}
from the asymptotic equality of 
\begin{equation} \label{ok}
\binom{n-1}{k-1}\sim \binom{n+\binom{k}{2}-1}{k-1},
\end{equation}
the leftmost and rightmost terms in  (\ref{guptaproof}), 
 under the assumption $k=o(n^{1/3})$.


To get the $n=k^3$ threshold function for integer partitions, we first show that for  ${n/k^3\rightarrow 0} $, 
\begin{eqnarray*}
k!q(n,k) &=& {\cal D}(n,k),\\ 
{\cal D}(n,k)&=&o( {\cal C}(n,k)               ), \hbox{ from (\ref{fenn}), }\\
\hbox { and }
{\cal C}(n,k)&\leq & k!p(n,k).
\end{eqnarray*}
To do the case 
$n/k^3\rightarrow \infty$, i.e. $k=o(n^{1/3})$, we  use Erd\H os--Lehner twice and also  (\ref{ok}):
$$q(n,k)=p\biggl(n-\binom{k}{2},k\biggl)\sim \frac{1}{k!}\binom{n-\binom{k}{2}-1}{k-1}\sim \frac{1}{k!}\binom{n-1}{k-1} \sim p(n,k).$$

\section{Threshold function for surjections} \label{sur}
Let $F(n,k)$ denote the number of $[n]\rightarrow [k]$ surjections. It is well-known  from the Bonferroni inequalities that
$k^n -n(k-1)^n  \leq F(n,k)\leq k^n$ and hence 
for ${n}>>{k\log k}$, $F(n,k)=
(1+o(1))k^n$ uniformly. Take an $f:[n]\rightarrow [k]$ random surjection. Let $A_{ij}(t)$ denote the event that for $1\leq i<j\leq k$ we have 
$|f^{-1}(i)|=|f^{-1}(j)|=t$. Observe that $P(A_{ij}(t))=\sum_{t\geq 1} \binom{n}{2t}\binom{2t}{t} \frac{F(n-2t,k-2)}{F(n,k)}$ and recall
$\binom{2t}{t}\sim \frac{2^{2t}}{\sqrt{\pi t}}$. Observe$^*$ 
\footnote{$^*$Note that for large $t$ (i.e. $n-2t= O(k \ln k)$) the approximation for  $P(A_{ij}(t))$ is not accurate. The same problem  occurs for small $t$  ($t=O(1)$), because of the estimate of $\binom{2t}{t}$. The corresponding terms, however, are negligible  both in the sum of probabilities and in (\ref{vmi}).}
that
\begin{eqnarray} \nonumber
P(\exists \hbox{$i<j$:  }|f^{-1}(i)|=|f^{-1}(j)| )& =& P\Biggl( \bigcup_{i<j}\bigcup_t A_{ij}(t) \Biggl) \leq \sum_{i<j}\sum_t P(A_{ij}(t))\sim
 {k\choose 2} \sum_{t\geq 1}\binom{n}{2t}\binom{2t}{t} \frac{(k-2)^{n-2t }}{ k^n }\\
&\sim&
{k\choose 2} \Biggl( 1- \frac{2}{k}  \Biggl)^{n}    \sum_{t\geq 1} \binom{n}{2t} \Biggl(  \frac{2}{k-2}  \Biggl)^{2t} \frac{1}{\sqrt{\pi t}}. \label{vmi}
\end{eqnarray}
Let $b(n,i) $ denote the term $\binom{n}{i} p^i(1-p)^{n-i}$ from the binomial distribution with $p=1-\frac{2}{k}$. It is easy to see that the core summation in (\ref{vmi}) is 
$$\sum_{t\geq 1} \binom{n}{2t} \Biggl(  \frac{2}{k-2}  \Biggl)^{2t} \frac{1}{\sqrt{\pi t}}=\Biggl( 1-\frac{2}{k}\Biggl)^{-n} \sum_{t\geq 1} \frac{b(n,2t)}{\sqrt{\pi t}}.$$
Observe that for this binomial distribution $\mu=np=\frac{2n}{k}$ and $\sigma < \sqrt{np}=\sqrt{\frac{2n}{k}}$.
Recall from \cite{as} the large deviation
inequality for sums of independent Bernoulli random variables:
$$P\Biggl(    |Y-\mu|>\epsilon \mu \Biggl)<2 e^{c_\epsilon \mu}, $$
where $c_\epsilon=\min \{ \ln (\epsilon^\epsilon (1+\epsilon)^{(1+\epsilon)}, \epsilon^2/2\}$. We select $\epsilon=\frac{1}{\ln n}$, with which for sufficiently large $n$,
$c_\epsilon=  \epsilon^2/2=\frac{1}{2\ln^2 n}$. Set $A=(1-\epsilon)\mu$ and $B=(1+\epsilon)\mu$. As $[A,B]$ includes the range where the normal convergence takes place,
$\sum_{A\leq t\leq B} b(n,t) \sim 1$. By Claim~\ref{elementary}, $\sum_{A\leq  2t\leq B} b(n,2t) \sim 1/2$. Also, if $A\leq t\leq B$, then ${t}\sim \frac{2n}{k}$. 
By the large deviation inequality above and $k<\sqrt{n}$ from (\ref{natural}), we obtain
$$\sum_{t\geq 1\atop t\notin [A,B]} b(n,t)=o\Biggl(\sqrt{\frac{k}{2n}}\sum_{A\leq t\leq B} b(n,t)\Biggl),$$
and combining with Claim~\ref{elementary} we obtain
$$\sum_{t\geq 1\atop t\notin [A,B]} b(n,t)=o\Biggl(\sqrt{\frac{k}{2n}}\sum_{A\leq 2t\leq B} b(n,2t)\Biggl).$$
Putting together these arguments:
$$\sum_{t\geq 1} \frac{b(n,2t)}{\sqrt{\pi t}}\sim \sqrt{\frac{k}{2\pi n}}\sum_{A\leq 2t\leq B} b(n,2t)\sim \sqrt{\frac{k}{2\pi n}}\sum_{t\geq 1} b(n,2t)
\sim \frac{1}{2}  \sqrt{\frac{k}{2\pi n}}\sum_{t\geq 1} b(n,t)\sim \frac{1}{2}  \sqrt{\frac{k}{2\pi n}}.
$$
\old{
-------------------------
Note that the core term   $\binom{n}{2t} \Bigl(  \frac{2}{k-2}  \Bigl)^{2t}$   in (\ref{vmi}) is just 
$\Bigl( 1+ \frac{2}{k-2}  \Bigl)^{-n} $
 times every second term in the binomial distribution
 $\binom{n}{t} \Bigl(  \frac{2}{k}  \Bigl)^{t}\Bigl( 1- \frac{2}{k}  \Bigl)^{n-t}$, which has expectation $\frac{2n}{k}$ and variance $<\sqrt{\frac{2n}{k}}$.
We have from the Binomial Theorem 
$\sum_{t} \binom{n}{t} \Bigl(  \frac{2}{k-2}  \Bigl)^{t}= 
\Bigl( 1+ \frac{2}{k-2}  \Bigl)^{n}$. 
Therefore, the summation in (\ref{vmi}) can be restricted for $t$'s satisfying $|t-\frac{2n}{k}|<\sqrt{\frac{2n}{k}}\log n$, without changing the asymptotic value of (\ref{vmi}),
because of the normal convergence of the binomial distribution. ????
Moreover, as consecutive terms in the binomial distribution 
are asymptotically equal to each other where their contribution is significant, 
$ \sum_{t\geq 1} \binom{n}{2t} \Bigl(  \frac{2}{k-2}  \Bigl)^{2t}\sim \frac{1+o(1)}{2} \Bigl( 1+ \frac{2}{k-2}  \Bigl)^{n}$, and by the normal convergence 
of the binomial distribution, it is  concentrated around $t=n/k$. 
}
We obtain  asymptotic formula for the upper bound with
\begin{equation} \label{num-func}
\sum_{i<j}\sum_t P(A_{ij}(t))=(1+o(1))\frac{k^2}{4} \sqrt{\frac{k}{2n\pi}} 
\end{equation}
which goes to zero 
as $n/k^5\rightarrow \infty$.

Next we use (\ref{also}) and (\ref{num-func}) to show that $P(\exists \hbox{$i<j$:  }|f^{-1}(i)|=|f^{-1}(j)| ) \rightarrow 1$
as $n/k^5\rightarrow 0$. The numerator of (\ref{also}) is the square of (\ref{num-func})
that grows to infinity. Therefore we can neglect the same term without square in the denominator of (\ref{also}).

A second negligible term
in the denominator arises if $i<j$ and $u<v$ make only 3 distinct indices (note that they must occur with the same $t$). The corresponding sum of the probabilities is estimated by
\begin{eqnarray} \nonumber
&& k{k-1\choose 2}\sum_{t\geq 1}  \binom{n}{3t} \frac{(3t)!}{(t!)^3}  \frac{F  (n-3t,k-3)}{   F (n,k)} 
 \leq   {k^3}\sum_{t\geq 1}    \binom{n}{3t}  \frac{(3t/e)^{3t} \sqrt{6\pi t}         }{(t/e)^{3t} (\sqrt{2\pi t})^3 }\frac{(k-3)^{n-3t}}{k^n} \\
& &\leq {k^3}\sum_{t\geq 1}    \binom{n}{3t} \Biggl( \frac{3}{k-3} \Biggl)^{3t} \frac{\sqrt{3}}{2\pi t} \frac{(k-3)^{n}}{k^n}.\label{only3-func}
\end{eqnarray}
We have from the Binomial Theorem 
$\sum_{t} \binom{n}{t} \Bigl(  \frac{3}{k-3}  \Bigl)^{t}
=
\Bigl( 1+ \frac{3}{k-3}  \Bigl)^{n}$.  Working with every third term in a binomial distribution with $p=\frac{3}{k}$ like we worked above with every second,
one obtains
\old{
Moreover, as consecutive terms in the binomial distribution 
are asymptotically equal to each other where their contribution is significant, 
$ \sum_{t\geq 1} \binom{n}{3t} \Bigl(  \frac{3}{k-3}  \Bigl)^{3t}\sim \frac{1+o(1)}{3} \Bigl( 1+ \frac{3}{k-3}  \Bigl)^{n}$, and by the normal convergence 
of the binomial distribution, it is  concentrated around $t=n/k$. We continue 
}
 the upper bound  for (\ref{only3-func})
$$
O\Biggl( \frac{k^4}{ n}\Biggl)=o\Biggl( \Bigl(\sum_{i<j}\sum_t P(A_{ij}(t))\Bigl)^2 \Biggl),
$$
using (\ref{num-func}) as $n/k^5\rightarrow 0$.


A third negligible term
in the denominator arises if $i<j$ and $u<v$  are 4 distinct indices, but the corresponding parts (set sizes) are all equal. The corresponding  term is
$$\binom{k}{4}
\sum_{t\geq 1}  \binom{n}{4t} \frac{(4t)!}{(t!)^4} \frac{  F  (n-4t,k-4)}{  F  (n,k)}. 
$$
This is easily estimated by $O(k^{11/2}/n^{3/2})$ like the estimate in
(\ref{only3-func}), and is similarly negligible compared to (\ref{num-func}) as  $n/k^5\rightarrow 0$. 
The significant term in the denominator is 
\begin{equation} \label{4-func}
 {k\choose 2}{k-2\choose 2}\sum_{\ell\geq 1}\sum_{t\geq 1,t\not=\ell}   \binom{n}{2\ell} \binom{2\ell}{\ell}\binom{n-2\ell}{2t} \binom{2t}{t} \frac{F  (n-2\ell-2t,k-4)}{   F (n,k)}
 \end{equation}
corresponding to the cases when for 4 distinct indices  2-2 parts (set sizes) are the same. 
The siginificant term will not change asymptotically
when we add the $t=\ell$ cases to the summation.

We do not repeat below arguments about the binomial distribution  which we went 
through before. So (\ref{4-func}) is asymptotically equal to
\begin{eqnarray} \label{mainterm-func}
&& \frac{k^4}{4}\sum_{\ell\geq 1} \binom{n}{2\ell} \frac{4^\ell}{\sqrt{\ell \pi}}\sum_{t\geq 1}\binom{n-2\ell}{2t} \frac{4^t}{\sqrt{t \pi}}
 \frac{(k-4)^{n-2\ell-2t}}{k^n}\\
 &\sim& \frac{k^4}{4}\sum_{\ell\geq 1} \binom{n}{2\ell} \frac{4^\ell}{\sqrt{\ell \pi}} \frac{(k-4)^{n-2\ell}}{k^n} \sum_{t\geq 1} \binom{n-2\ell}{2t}  \Biggl( 
 \frac{2}{k-4} \Biggl)^{2t} \frac{1}{\sqrt{t \pi}}\\
 &\sim& \frac{k^4}{4}\sum_{\ell\geq 1} \binom{n}{2\ell} \frac{4^\ell}{\sqrt{\ell \pi}} \frac{(k-4)^{n-2\ell}}{k^n}\cdot \frac{1}{2} \Biggl( 
 1+\frac{2}{k-4} \Biggl)^{n-2\ell} \sqrt{\frac{k-2}{2\pi (n-2\ell)}}\\
 &\sim& \frac{k^4}{8\pi }  \Biggl( 
 1-\frac{2}{k} \Biggl)^{n} \sum_{\ell\geq 1} \binom{n}{2\ell} \Biggl( 
 \frac{2}{k-2} \Biggl)^{2\ell} \sqrt{\frac{k-2}{2\ell (n-2\ell)}}\\
 &\sim& \frac{k^4}{8\pi }  \Biggl( 
 1-\frac{2}{k} \Biggl)^{n} \cdot \frac{1}{2} \Biggl( 
 1+\frac{2}{k-2} \Biggl)^{n}\sqrt{\frac{k-2}{\frac{4n}{k} (n-4\frac{n}{k})}}\sim \frac{k^5}{32\pi n}
  \end{eqnarray}
  We conclude that  the  numerator in (\ref{also}), which is (\ref{num-func}) squared in our setting, is asymptotically equal to its denominator (\ref{mainterm-func}), proving that
  $P(\exists \hbox{$i<j$:  }|f^{-1}(i)|=|f^{-1}(j)| )\rightarrow 1$
as $n/k^5\rightarrow 0$.

\end{document}